\documentclass[12pt,a4paper]{article}
\usepackage[UTF8]{ctex}
\usepackage{amsmath,amssymb,amsthm,mathrsfs}
\usepackage{geometry,graphicx,booktabs,longtable}
\usepackage[hidelinks]{hyperref}
\allowdisplaybreaks[3]
\newtheorem{theorem}{Theorem}[section]
\newtheorem{proposition}[theorem]{Proposition}
\newtheorem{lemma}[theorem]{Lemma}

\theoremstyle{definition}

\newtheorem{definition}[theorem]{Definition}
\theoremstyle{plain}
\newtheorem{extheorem}{Theorem}

\numberwithin{equation}{section}
\DeclareMathOperator{\Ent}{Ent}
\DeclareMathOperator{\Rea}{Re}
\DeclareMathOperator{\Ima}{Im}
\DeclareMathOperator{\sech}{sech}
\newcommand{\C}{\mathbb C}
\newcommand{\R}{\mathbb R}
\newcommand{\D}{\overline{\mathbb D}}
\newcommand{\norm}[1]{\lVert#1\rVert}
\newcommand{\abs}[1]{\lvert#1\rvert}
\newcommand{\E}{\mathbb E_\nu}
\graphicspath{{figures/}}

\begin{document}
\title{The Two-Point Inequality for $p>q$}
\author{Di Wu\thanks{School of Science, Zhejiang University of Science and Technology,
Hangzhou 310023, P.R. China. E-mail address: \texttt{wudi08@mails.ucas.ac.cn}.
Supported by the Natural Science Foundation of Zhejiang Province of China (Grant
No.~LQ18A010002) and the National Natural Science Foundation of China (Grant
No.~12001488).}}
\date{}
\maketitle
\begin{abstract}
For $1<q<p<\infty$ we study the $L^p\to L^q$ contractivity of the two-point operator
$T_\omega(a+bx)=a+\omega bx$ with a complex parameter $\omega$,
extending the Weissler two-point condition, usually discussed in the range $1<p\le q<\infty$,
to the reversed order of exponents.
For $1<q<p<2$ we prove that the second-order necessary condition $\abs{(p-2)-(q-2)\omega^2}\le p-q\abs{\omega}^2$,
obtained near $(a,b)=(1,0)$, is also sufficient;
by duality the same description applies to $2<q<p<\infty$.
The proof uses the method of characteristic flows. In addition,
for $q\le2\le p$ the contractivity follows directly from $L^2$ contraction.
\end{abstract}


\section{Introduction}\label{sec:intro}

On the two-point space $\{-1,1\}$ we take the probability measure
$$
d\nu=\frac{\delta_{-1}+\delta_1}{2}.
$$
Every complex-valued function is uniquely written as $f(x)=a+bx$ with $a,b\in\C$.
For $1<r<\infty$ and $\omega\in\C$ define
$$
\norm{f}_r=\left(\frac{\abs{a+b}^r+\abs{a-b}^r}{2}\right)^{1/r},
\qquad T_\omega f(x)=a+\omega bx.
$$
We study the range of complex parameters for which the two-point inequality
\begin{equation}\label{eq:twopoint}
\left(\frac{\abs{a+\omega b}^q+\abs{a-\omega b}^q}{2}\right)^{1/q}
\le
\left(\frac{\abs{a+b}^p+\abs{a-b}^p}{2}\right)^{1/p}
\quad(a,b\in\C)
\end{equation}
holds. Write
\begin{align}
\Omega_{p,q}
&=\{\omega\in\C:\norm{T_\omega f}_q\le\norm{f}_p\text{ for all }f\},\label{eq:Omega}\\
A_{p,q}
&=\{\omega\in\C:\abs{(p-2)-(q-2)\omega^2}\le p-q\abs{\omega}^2\}.\label{eq:A}
\end{align}
For $(a,b)=(1,0)$ both sides of the two-point inequality are equal to $1$,
so for any choice of the parameters the operator norm is at least $1$.
Consequently $\omega\in\Omega_{p,q}$ means precisely that the best constant equals $1$.

In the range $1<p\le q<\infty$ the two-point inequality is closely related to the Hermite semigroup
and to sharp Fourier inequalities.
Beckner's discrete two-point method was an important step in obtaining the sharp Hausdorff--Young
inequality \cite{Beckner1975}.
In this range of exponents Weissler proved the sufficiency of condition~\eqref{eq:A},
but left open the two ranges $2<p\le q<3$ and $3/2<p\le q<2$ \cite{Weissler1979}.
For $p\le q$ two-point contraction can be transferred by tensorization to Hamming cubes of arbitrary dimension.
Ivanisvili--Nazarov proved the case $1<p=q<\infty$
\cite{IvanisviliNazarov2022}.
Huang--Ivanisvili gave, in the preprint \cite{HuangIvanisvili2026},
a proof by the flow method for $2<p\le q<\infty$, and obtained $1<p\le q<2$ by duality.
The theorem stated in that paper covers the two ranges of exponents left open above,
and the case $p=q$ is also handled directly by their method.

The present paper is concerned with the reversed order of exponents $p>q$.
For $1<q<p<2$ or $2<q<p<\infty$,
a second-order variation near $(a,b)=(1,0)$ again yields condition~\eqref{eq:A}.
We consider contractivity only on the two-point space; for the present range $p>q$,
we make no tensorization claim analogous to that in the range $p\le q$.
Our main result is the following.

\begin{theorem}\label{thm:main}
Let either $1<q<p<2$ or $2<q<p<\infty$. Then
$$
\Omega_{p,q}=A_{p,q}.
$$
In other words, \eqref{eq:twopoint} holds for all $a,b\in\C$ if and only if
$$
\abs{(p-2)-(q-2)\omega^2}\le p-q\abs{\omega}^2.
$$
For every $\omega\in A_{p,q}$, the operator norm of $T_\omega:L^p(d\nu)\to L^q(d\nu)$ is $1$,
attained at $(a,b)=(1,0)$.
\end{theorem}

The proof of this theorem uses the method of characteristic flows. More precisely,
we fix the input exponent $p\in(1,2)$, let the output exponent $s$ decrease from $p$,
and move the complex parameter $\omega$ so that $\omega$ always lies on $\partial A_{p,s}$ and
$\norm{T_\omega f}_s$ is non-increasing in forward time for every two-point function $f$.
Solving this flow backwards from the target boundary and then taking a subsequential limit as $s\uparrow p$,
we can apply the following known result.

\begin{extheorem}[Ivanisvili--Nazarov {\cite{IvanisviliNazarov2022}}]\label{thm:diagonal}
For all $1<p=q<\infty$ one has $\Omega_{p,p}=A_{p,p}$.
\end{extheorem}

Before describing the local mechanism of the proof we recall the standard notation for entropy.
For a non-negative two-point function $h$ set
$$
\Ent_\nu(h)=\E(h\log h)-(\E h)\log(\E h),
\qquad 0\log0=0.
$$
Let $P(a+bx)=bx$. Since $a=\E f$, we have $Pf=f-\E f$, that is, $P$ removes the mean of $f$.
\begin{definition}[Admissible set of complex directions]\label{def:control}
For $1<s<2$ set
\begin{equation}\label{eq:control}
\mathcal K_s=
\left\{-\alpha+i\beta:\alpha\ge0,\quad
(2-s)^2\beta^2\le2\alpha+4(s-1)\alpha^2\right\}.
\end{equation}
\end{definition}

\begin{lemma}[Complex two-point inequality of logarithmic Sobolev type]\label{lem:ITP}
Let $1<s<2$, $\eta\in\mathcal K_s$, and let $f$ be any complex-valued two-point function. Then
\begin{equation}\label{eq:LSI}
\Rea\!\left(\eta\,\E\!\left[\abs{f}^{s-2}\overline f\,Pf\right]\right)
\le\frac1{s^2}\Ent_\nu(\abs{f}^s).
\end{equation}
At points where $f=0$ the expression $\abs{f}^{s-2}\overline f$ is defined to be $0$ by continuity.
Equivalently, for every non-zero $f$,
\begin{equation}\label{eq:ITP}
\left.\frac{d}{d\varepsilon}\right|_{\varepsilon=0}
\log\norm{T_{1+\varepsilon\eta}f}_{s-\varepsilon}\le0.
\end{equation}
\end{lemma}

\begin{lemma}[A direction tangent to the boundary surface]\label{lem:tangent}
For $p,s>1$ and $\omega\in\C$ write
\begin{equation}\label{eq:F}
\Lambda(p,s,\omega)=(p-2)-(s-2)\omega^2,\qquad
F(p,s,\omega)=p-s\abs{\omega}^2-\abs{\Lambda(p,s,\omega)}.
\end{equation}
When $\omega\ne0$ and $\Lambda\ne0$ define
\begin{equation}\label{eq:csigma}
c+i\sigma=\frac{\overline\Lambda\,\omega^2}{\abs{\Lambda}\abs{\omega}^2},
\qquad c^2+\sigma^2=1,
\end{equation}
and the first-order change of the boundary function
\begin{equation}\label{eq:Ddef}
\mathcal D_{p,s,\omega}(\eta)=
\left.\frac{d}{d\varepsilon}\right|_{\varepsilon=0}
F\bigl(p,s-\varepsilon,\omega(1+\varepsilon\eta)\bigr).
\end{equation}

If $1<s<p<2$ and $F(p,s,\omega)=0$, then $\omega\ne0$, $\Lambda\ne0$, and
$$
sc+2-s>0.
$$
On the admissible set $\mathcal K_s$ the functional $\mathcal D_{p,s,\omega}$ has minimum $0$;
the minimum is attained in the direction
\begin{equation}\label{eq:eta}
\eta_*(\omega,s)=
-\frac{1-c}{2(sc+2-s)}
-i\,\frac{\sigma}{(2-s)(sc+2-s)}\in\partial\mathcal K_s.
\end{equation}
On the open set where $1<s<2$, $\omega\ne0$,
$\Lambda\ne0$ and $sc+2-s>0$ this formula is smooth in the real variables $(\Rea\omega,\Ima\omega,s)$.
\end{lemma}

These two lemmas deal with the norm and with the boundary of the parameter region respectively.
Lemma~\ref{lem:ITP} says that, along an admissible direction,
changing the multiplier and the exponent simultaneously makes the first-order change of the norm non-positive;
Lemma~\ref{lem:tangent} selects, among these directions,
a direction for which the first-order change of the boundary function vanishes. Hence the system
\begin{equation}\label{eq:system}
\omega'(t)=\omega(t)\eta_*(\omega(t),s(t)),\qquad s'(t)=-1
\end{equation}
satisfies both
$$
\frac{d}{dt}F(p,s(t),\omega(t))=0,\qquad
\frac{d}{dt}\norm{T_{\omega(t)}f}_{s(t)}\le0.
$$
The first identity keeps the trajectory on the boundary surface;
the second provides the monotonicity needed to propagate contractivity from $s=p$ to $1<s<p<2$.
The starting result at $s=p$ is supplied by Theorem~\ref{thm:diagonal}.
The flow is chosen so as to meet these two requirements,
namely preservation of the boundary and monotonicity of the norm.
It is a characteristic flow in parameter space;
different functions $f$ share the same parameter orbit.
We take the contraction statement for $p=q$ as an input;
starting merely from $A_{p,2}$ while keeping monotonicity of the norm for every $f$ does not cover the whole target set,
so this input cannot be removed in this way. On the other hand,
for $p<s<2$ the regularity condition $sc+2-s>0$ on the boundary no longer holds everywhere,
so the present construction cannot be applied directly to the range $p<q$ and does not directly yield the result of
\cite{HuangIvanisvili2026}. These two limitations are discussed in
Section~\ref{sec:scope}. Another manuscript \cite{Wu2026}
gives a proof which does not rely on characteristic flows and treats the full range $1\le p,q\le\infty$.

\section{Preliminaries and proofs of the main lemmas}\label{sec:lemmas}

\subsection{The second-order condition, the unit disc and duality}

For $1<r<\infty$ and $v\in\C$ set
$$
Q_r(v)=(r-1)(\Rea v)^2+(\Ima v)^2.
$$
Expanding near $(a,b)=(1,0)$ gives
\begin{equation}\label{eq:Taylor}
\norm{1+\varepsilon vx}_r^2=1+\varepsilon^2Q_r(v)+o(\varepsilon^2),
\qquad \varepsilon\in\R,\quad\varepsilon\to0.
\end{equation}
Indeed,
if we expand $\abs{1+\varepsilon v}^r$ and $\abs{1-\varepsilon v}^r$ and average the two expansions,
the odd-order terms cancel and the coefficient of $\varepsilon^2$ is $rQ_r(v)/2$.
Raising the result to the power $2/r$ then yields the formula.

\begin{proposition}\label{prop:local}
For $1<p,q<\infty$,
\begin{equation}\label{eq:Q}
\omega\in A_{p,q}
\quad\Longleftrightarrow\quad
Q_q(\omega v)\le Q_p(v)\quad(v\in\C).
\end{equation}
Consequently $\Omega_{p,q}\subseteq A_{p,q}$. If $1<q\le p<2$, then $A_{p,q}\subseteq\D$.
For fixed $p$ and $\omega$, the function $s\mapsto F(p,s,\omega)$ is non-increasing.
\end{proposition}

\begin{proof}
Write $\omega=u+iv$ and let $W$ be the matrix of the quadratic form $Q_p(z)-Q_q(\omega z)$; then
$$
W=\begin{pmatrix}
p-1-(q-1)u^2-v^2&(q-2)uv\\
(q-2)uv&1-u^2-(q-1)v^2
\end{pmatrix}.
$$
Its trace and the square of the difference of its eigenvalues satisfy
$$
\operatorname{tr}W=p-q\abs{\omega}^2,\qquad
(W_{11}-W_{22})^2+4W_{12}^2=\abs{\Lambda(p,q,\omega)}^2.
$$
Hence the smaller eigenvalue equals exactly $F(p,q,\omega)/2$, which proves
\eqref{eq:Q}.
Substituting \eqref{eq:Taylor} into the contraction inequality gives the necessity.

If $q<2$, then by the reverse triangle inequality and $F\ge0$,
$$
p-q\abs{\omega}^2
\ge\abs{(2-q)\omega^2-(2-p)}
\ge(2-q)\abs{\omega}^2-(2-p),
$$
and rearranging yields $\abs{\omega}^2\le1$. Finally, if $s_2\ge s_1$, then
$$
\abs{\Lambda(p,s_2,\omega)}-\abs{\Lambda(p,s_1,\omega)}
\ge-(s_2-s_1)\abs{\omega}^2.
$$
Substituting this into $F(p,s_2,\omega)-F(p,s_1,\omega)$ shows that the difference is at most zero.
The argument also covers the points where $\Lambda=0$.
\end{proof}

Equivalence~\eqref{eq:Q} also shows that $A_{p,q}$ is a compact convex set containing the origin,
which is an interior point:
for every non-zero $v$ the condition $Q_q(\omega v)\le Q_p(v)$ defines a closed elliptical disc,
and positive definiteness guarantees that it holds strictly near the origin.
Along any ray $\omega=re^{i\theta}$ this condition takes the form $r^2Q_q(e^{i\theta}v)\le Q_p(v)$;
hence there is exactly one boundary radius on each ray, and
\begin{equation}\label{eq:boundary}
\partial A_{p,q}=\{\omega:F(p,q,\omega)=0\}.
\end{equation}

\begin{lemma}[Duality]\label{lem:duality}
If $p'=p/(p-1)$ and $q'=q/(q-1)$, then
$$
\Omega_{p,q}=\Omega_{q',p'},\qquad A_{p,q}=A_{q',p'}.
$$
\end{lemma}
\begin{proof}
For the pairing $\langle f,g\rangle=\E(f\overline g)$ we have $T_\omega^*=T_{\overline\omega}$,
hence $\norm{T_\omega}_{p\to q}=\norm{T_{\overline\omega}}_{q'\to p'}$.
Complex conjugation preserves all two-point norms, which yields the identity for $\Omega$.

To prove the identity for $A$,
let $D_r=\operatorname{diag}(r-1,1)$ and let $L_\omega$ denote the real matrix of multiplication by $\omega$.
Condition~\eqref{eq:Q} is equivalent to
$$
\norm{D_q^{1/2}L_\omega D_p^{-1/2}}_{2\to2}\le1.
$$
Transposition does not change the Euclidean operator norm,
and $L_\omega^{\mathsf T}=L_{\overline\omega}$, $D_{r'}=D_r^{-1}$.
The transposed condition is exactly $\overline\omega\in A_{q',p'}$.
Using the conjugation symmetry once more gives the claim.
\end{proof}

\subsection{Derivatives of two-point norms and entropy}

\begin{lemma}[Derivative formula]\label{lem:derivative}
Let $s>1$, let $f$ be a fixed non-zero two-point function and let $\eta\in\C$. Define
\begin{equation}\label{eq:EV}
E_s(f)=\frac{\Ent_\nu(\abs{f}^s)}{s^2\norm{f}_s^s},
\qquad
V_s(f)=\frac{\E\left[\abs{f}^{s-2}\overline f\,Pf\right]}{\norm{f}_s^s}.
\end{equation}
Then
\begin{equation}\label{eq:derivative}
\left.\frac{d}{d\varepsilon}\right|_{\varepsilon=0}
\log\norm{T_{1+\varepsilon\eta}f}_{s-\varepsilon}
=\Rea(\eta V_s(f))-E_s(f),
\qquad
\partial_s\log\norm{f}_s=E_s(f).
\end{equation}
\end{lemma}
\begin{proof}
First suppose that $f$ does not vanish at either point and write $M=\E\abs{f}^s$.
Since $T_{1+\varepsilon\eta}f=f+\varepsilon\eta Pf$,
differentiating $(s-\varepsilon)^{-1}\log\E\abs{f+\varepsilon\eta Pf}^{s-\varepsilon}$ directly gives
\begin{align*}
\left.\frac{d}{d\varepsilon}\right|_0
\log\norm{T_{1+\varepsilon\eta}f}_{s-\varepsilon}
&=\frac{\log M}{s^2}
-\frac{\E(\abs{f}^s\log\abs{f})}{sM}
+\frac{\Rea\!\left(\eta\,\E[\abs{f}^{s-2}\overline f\,Pf]\right)}{M}\\
&=-\frac{\Ent_\nu(\abs{f}^s)}{s^2M}+\Rea(\eta V_s(f)).
\end{align*}
Differentiating $s^{-1}\log\E\abs{f}^s$ in the same way gives the second identity.

If $f$ vanishes at one point,
the varying term at that point is of the form $\abs{\varepsilon z}^{s-\varepsilon}=o(\abs\varepsilon)$;
taking $\varepsilon$ so small that $s-\varepsilon>1$ gives a two-sided derivative equal to zero there.
Moreover $u^s\log u\to0$ as $u\downarrow0$,
so the convention for the entropy term is consistent as well. More generally,
$(r,z)\mapsto\abs z^r$ is jointly $C^1$ on $(1,\infty)\times\C$,
so the formulas above hold for every non-zero $f$.
\end{proof}

By the derivative formula, \eqref{eq:ITP} is equivalent to
$$
\Rea(\eta V_s(f))\le E_s(f).
$$
Multiplying both sides by $\norm f_s^s>0$ gives exactly \eqref{eq:LSI}. Conversely,
dividing by this positive number and substituting
\eqref{eq:derivative} recovers \eqref{eq:ITP}; for $f=0$,
\eqref{eq:LSI} holds separately. This shows that the two formulations are
equivalent.

We call this a complex two-point inequality of logarithmic Sobolev type because the variation of the
exponent produces the standard entropy term,
whereas the variation of the multiplier produces the generator pairing term $\E[\abs f^{s-2}\overline f\,Pf]$.
Indeed, $T_{e^{-u}}=e^{-uP}$,
so $P$ is the non-negative generator of the corresponding semigroup $(e^{-uP})_{u\ge0}$.
The complex parameter $\eta$ here couples a dissipative direction with a change of phase,
and the sense in which the exponent varies is fixed by $s'=-1$; hence \eqref{eq:LSI}
is an entropy estimate adapted to this variable-exponent flow.
The direction of its inequality should be distinguished from that in the usual real-valued $L^2$
logarithmic Sobolev inequality.
The non-negativity of the entropy itself comes from the convexity of $u\mapsto u\log u$.

For $f=1+\zeta x$ write
$$
u_\pm=\abs{1\pm\zeta},\qquad
w_\pm=\frac{u_\pm^s}{u_+^s+u_-^s}.
$$
Then
\begin{equation}\label{eq:weights}
E_s(f)=\frac{w_+\log(2w_+)+w_-\log(2w_-)}{s^2},
\qquad
V_s(f)=w_+\frac{\zeta}{1+\zeta}-w_-\frac{\zeta}{1-\zeta}.
\end{equation}
If a denominator vanishes, the corresponding weighted term as a whole is defined to be zero.
For instance, for $\zeta=1$ we have $V_s(f)=1/2$.

\subsection{Proof of the complex two-point inequality of logarithmic Sobolev type}

\begin{proof}[Proof of Lemma~\ref{lem:ITP}]
If $f=0$, \eqref{eq:LSI} is trivial. If $f=bx$ with $b\ne0$,
then $\norm{T_{1+\varepsilon\eta}f}_{s-\varepsilon}=\abs{1+\varepsilon\eta}\abs b$,
whose logarithmic derivative is $\Rea\eta=-\alpha\le0$.
In the remaining cases homogeneity allows us to assume $f=1+\zeta x$. For $\zeta=\pm1$,
\eqref{eq:weights} gives $V_s=1/2$ and $E_s=\log2/s^2$, and the claim follows.
Henceforth assume $\zeta\ne\pm1$.

\medskip\noindent\textit{Coordinates given by the ratio of the two values.}
To simplify both the weights determined by the moduli and the phase, write
\begin{equation}\label{eq:bipolar}
\frac{1+\zeta}{1-\zeta}=e^{L+i\Theta},
\qquad L\in\R,\quad \Theta\in(-\pi,\pi].
\end{equation}
Then $\zeta=\tanh((L+i\Theta)/2)$ and
$$
y=\frac{sL}{2},\qquad h=\frac{(2-s)L}{2},\qquad
w_\pm=\frac{e^{\pm y}}{2\cosh y},\qquad w_+-w_-=\tanh y.
$$
Here $L$ is the logarithm of the ratio of the moduli of the two values and $\Theta$ is their phase difference.
Define
\begin{equation}\label{eq:Phi}
\Phi(y)=y\tanh y-\log\cosh y.
\end{equation}
Since $2w_\pm=e^{\pm y}/\cosh y$,
$$
E_s(f)=\frac{\Phi(y)}{s^2}.
$$
To compute $V_s$, put \eqref{eq:weights} over a common denominator and use
$$
\frac{\zeta}{1-\zeta^2}=\frac12\sinh(L+i\Theta),\qquad
\frac{\zeta^2}{1-\zeta^2}=\frac{\cosh(L+i\Theta)-1}{2},
$$
which gives
$$
V_s=\frac12\bigl[1+\tanh y\,\sinh(L+i\Theta)-\cosh(L+i\Theta)\bigr].
$$
The hyperbolic subtraction formulas further give
\begin{align*}
\cosh L-\tanh y\,\sinh L&=\frac{\cosh(L-y)}{\cosh y}
=\frac{\cosh h}{\cosh y},\\
\tanh y\,\cosh L-\sinh L&=-\frac{\sinh(L-y)}{\cosh y}
=-\frac{\sinh h}{\cosh y}.
\end{align*}
Therefore
\begin{equation}\label{eq:coordinates}
\Rea V_s=\frac12\left(1-\cos\Theta\frac{\cosh h}{\cosh y}\right),
\qquad
\Ima V_s=-\frac{\sin\Theta\,\sinh h}{2\cosh y},
\qquad E_s=\frac{\Phi(y)}{s^2}.
\end{equation}
Since $\Phi(0)=0$ and $\Phi'(y)=y\sech^2y$, the function $\Phi$ is non-negative and even.
The quantities $\Rea V_s$, $\abs{\Ima V_s}$ and $E_s$ in
\eqref{eq:coordinates} are invariant under $L\mapsto-L$, so in the sequel we
may assume $L\ge0$ and hence $0\le h\le y$. In particular $\Rea V_s\ge0$.

\medskip\noindent\textit{Maximizing over the control direction and the phase.}
Fix $\alpha\ge0$. By \eqref{eq:control},
$$
\Rea(\eta V_s)=-\alpha\Rea V_s-\beta\Ima V_s
\le-\alpha\Rea V_s+
\frac{\abs{\Ima V_s}}{2-s}\sqrt{2\alpha+4(s-1)\alpha^2}.
$$
Thus it suffices to prove
\begin{equation}\label{eq:worstbeta}
\abs{\Ima V_s}\sqrt{2\alpha+4(s-1)\alpha^2}
\le(2-s)(E_s+\alpha\Rea V_s).
\end{equation}
To organize the quadratic expressions below, set
$$
B_s(y)=2\cosh y\,E_s=\frac{2}{s^2}\cosh y\,\Phi(y)\ge0.
$$
We substitute \eqref{eq:coordinates} into \eqref{eq:worstbeta}, multiply by $2\cosh y>0$,
and rearrange to obtain the equivalent inequality
\begin{align}
&\abs{\sin\Theta}\sinh h\sqrt{2\alpha+4(s-1)\alpha^2}
+(2-s)\alpha\cosh h\cos\Theta \notag\\
&\hspace{35mm}\le(2-s)\bigl(B_s(y)+\alpha\cosh y\bigr).\label{eq:phase}
\end{align}
By the Cauchy--Schwarz inequality and $\sin^2\Theta+\cos^2\Theta=1$, the left-hand side is at most
$$
\sqrt{\sinh^2h\,[2\alpha+4(s-1)\alpha^2]
+(2-s)^2\alpha^2\cosh^2h}.
$$
The right-hand side of \eqref{eq:phase} is non-negative, so we may compare squares.
The square of the right-hand side minus the quantity under the square root equals exactly
$$
c_2\alpha^2+c_1\alpha+c_0,
$$
where
\begin{align}
c_2&=(2-s)^2\sinh^2y-s^2\sinh^2h,\label{eq:c2}\\
c_1&=2\left[\frac{2(2-s)^2}{s^2}\cosh^2y\,\Phi(y)-\sinh^2h\right],\label{eq:c1}\\
c_0&=(2-s)^2B_s(y)^2.\label{eq:c0}
\end{align}
In collecting the terms of order $\alpha^2$ we used $(2-s)^2+4(s-1)=s^2$ and $\cosh^2u-1=\sinh^2u$.

\medskip\noindent\textit{Non-negativity of the three coefficients.}
Let $\kappa=(2-s)/s\in(0,1)$; then $h=\kappa y$.
Since $\sinh$ is convex on $[0,\infty)$ and $\sinh0=0$,
\begin{equation}\label{eq:sinhconvex}
\sinh(\kappa y)\le \kappa\sinh y.
\end{equation}
This immediately gives $c_2\ge0$. To treat $c_1$ we first prove
\begin{equation}\label{eq:entropybound}
2\Phi(y)\ge\tanh^2y\qquad(y\ge0).
\end{equation}
The difference of the two sides vanishes at the origin, and its derivative is
$$
\frac{d}{dy}\bigl(2\Phi(y)-\tanh^2y\bigr)
=2\sech^2y\,(y-\tanh y)\ge0.
$$
The last inequality follows from $(y-\tanh y)'=\tanh^2y\ge0$. Therefore
$$
\frac{c_1}{2}
=2\kappa^2\cosh^2y\,\Phi(y)-\sinh^2(\kappa y)
\ge \kappa^2\sinh^2y-\sinh^2(\kappa y)\ge0.
$$
Also $c_0\ge0$, so the required quadratic expression is non-negative for all $\alpha\ge0$.
Substituting back successively gives \eqref{eq:phase}, \eqref{eq:worstbeta}
and $\Rea(\eta V_s)\le E_s$. Combining this with Lemma~\ref{lem:derivative} completes the proof.
\end{proof}

The admissible set \eqref{eq:control}
can also be read off from a linearization of the second-order condition.
Substituting $1+\varepsilon\eta$ into $A_{s,s-\varepsilon}$,
where $\eta=-\alpha+i\beta$ and $\varepsilon\downarrow0$, gives
\begin{align*}
s-(s-\varepsilon)\abs{1+\varepsilon\eta}^2
&=\varepsilon(1+2s\alpha)+O(\varepsilon^2),\\
(s-2)-(s-\varepsilon-2)(1+\varepsilon\eta)^2
&=\varepsilon\bigl(1+2(2-s)\eta\bigr)+O(\varepsilon^2).
\end{align*}
Since $A_{s,s-\varepsilon}\subseteq\D$,
a necessary condition for $1+\varepsilon\eta$ to be admissible is $\abs{1+\varepsilon\eta}\le1$.
Its first-order expansion as $\varepsilon\downarrow0$ gives $\alpha\ge0$. Under this assumption,
comparing the first-order terms and squaring gives exactly $(2-s)^2\beta^2\le2\alpha+4(s-1)\alpha^2$.
This is the origin of $\mathcal K_s$; the proof of Lemma~\ref{lem:ITP}
verifies independently that this first-order condition indeed controls all two-point functions.

\subsection{Boundary-tangent direction: computation and regularity}

We now compute $\mathcal D_{p,s,\omega}$. In the computation below we assume only $1<s<2$,
$\omega\ne0$ and $\Lambda\ne0$; $F=0$ is not required. By \eqref{eq:F},
$$
-\partial_sF=\abs{\omega}^2-
\frac{\Rea(\overline\Lambda\,\omega^2)}{\abs\Lambda}
=\abs{\omega}^2(1-c).
$$
On the other hand,
under $\omega\mapsto\omega(1+\varepsilon\eta)$ the first-order change of $\Lambda$ equals $2(2-s)\omega^2\eta$.
Writing $\eta=-\alpha+i\beta$ therefore gives
\begin{equation}\label{eq:D}
\mathcal D_{p,s,\omega}(-\alpha+i\beta)
=\abs{\omega}^2\bigl[(1-c)+2(s+(2-s)c)\alpha
+2(2-s)\sigma\beta\bigr].
\end{equation}

\begin{proof}[Proof of Lemma~\ref{lem:tangent}]
\textit{Minimizing the first-order change.}
Assume for the moment that $sc+2-s>0$ and write $A=s+(2-s)c$. Then $A\ge2s-2>0$. For fixed $\alpha$,
the minimum of \eqref{eq:D} over the admissible $\beta$ is
$$
\abs{\omega}^2
\left[(1-c)+2A\alpha-2\abs\sigma
\sqrt{2\alpha+4(s-1)\alpha^2}\right].
$$
To prove that this expression is non-negative, consider
\begin{align}
R(\alpha)
&=\left(A\alpha+\frac{1-c}{2}\right)^2
-\sigma^2\bigl(2\alpha+4(s-1)\alpha^2\bigr)\notag\\
&=\bigl(A^2-4(s-1)\sigma^2\bigr)\alpha^2
+\bigl(A(1-c)-2\sigma^2\bigr)\alpha+\frac{(1-c)^2}{4}.\label{eq:poly}
\end{align}
Its discriminant is
\begin{align*}
\Delta
&=\bigl(A(1-c)-2\sigma^2\bigr)^2
-\bigl(A^2-4(s-1)\sigma^2\bigr)(1-c)^2\\
&=\sigma^2(1-c)\bigl[-4A+4(1+c)+4(s-1)(1-c)\bigr]=0.
\end{align*}
In the last step we substituted $A=s+(2-s)c$;
the constant terms and the coefficients of $c$ cancel separately.

To obtain the minimizing direction and the smoothness at the same time,
we use $\sigma^2=1-c^2$ once more and rearrange:
\begin{align}
A^2-4(s-1)\sigma^2&=(sc+2-s)^2,\notag\\
2\sigma^2-A(1-c)&=(1-c)(sc+2-s).\label{eq:simplify}
\end{align}
Thus \eqref{eq:poly} is in fact a perfect square,
$$
R(\alpha)=\left((sc+2-s)\alpha-\frac{1-c}{2}\right)^2.
$$
Since $A\alpha+(1-c)/2\ge0$, comparing squares gives back
$$
A\alpha+\frac{1-c}{2}\ge
\abs\sigma\sqrt{2\alpha+4(s-1)\alpha^2}.
$$
Equality holds for
$$
\alpha_*=\frac{1-c}{2(sc+2-s)}\ge0.
$$
Substituting this directly and using $\sigma^2=1-c^2$ yields
$$
2\alpha_*+4(s-1)\alpha_*^2
=\frac{\sigma^2}{(sc+2-s)^2}.
$$
Hence, taking $\beta_*=-\sigma/[(2-s)(sc+2-s)]$ produces \eqref{eq:eta},
which belongs to $\partial\mathcal K_s$ and makes $\mathcal D=0$.
At the same time $\mathcal D(\eta)\ge0$ for every $\eta\in\mathcal K_s$.
This conclusion holds on the whole of the open set above, and $F=0$ has not been used.

\medskip\noindent\textit{Boundary points satisfy the regularity condition.}
Under $1<s<p<2$,
$$
F(p,s,0)=2p-2>0,
$$
so $\omega\ne0$ on the boundary. If $\Lambda=0$, then $\omega^2=(2-p)/(2-s)$, and hence
$$
F(p,s,\omega)=p-s\frac{2-p}{2-s}
=\frac{2(p-s)}{2-s}>0,
$$
so such a point cannot lie on the boundary either. Write $\omega=\abs\omega e^{i\theta}$.
From $F=0$ we get $\abs\Lambda=p-s\abs\omega^2$. Therefore
\begin{align*}
\abs\Lambda\abs\omega^2(sc+2-s)
&=s\Rea(\overline\Lambda\,\omega^2)+(2-s)\abs\Lambda\abs\omega^2\\
&=\abs\omega^2\bigl[(2-s)p-s(2-p)\cos2\theta\bigr]\\
&\ge2\abs\omega^2(p-s)>0.
\end{align*}
Since $\abs\Lambda\abs\omega^2$ is strictly positive, we obtain $sc+2-s>0$.
This is precisely where the order of the exponents $p>s$ is used.

Finally, $c$ and $\sigma$ are smooth in the real variables on the set $\omega\ne0$, $\Lambda\ne0$,
and the denominators in \eqref{eq:eta} do not vanish on the open set described.
This also covers the points where $\sigma=0$,
so no separate branch for the sign function or the square root is needed.
\end{proof}

The term \emph{extremal direction} has a precise meaning here:
$\eta_*$ is the direction in the admissible set for which the affine functional $\mathcal D$ attains its minimum $0$.
The corresponding velocity $\omega\eta_*$, together with the exponent velocity $-1$,
gives a tangential direction of the boundary surface.
Its trajectories in the complex plane and the direction arrows are shown in Figure~\ref{fig:flow}(a);
the planar velocity is in general not perpendicular to the boundary at a fixed exponent,
see Section~\ref{sec:geometry}.

\subsection{Continuation of the ODE}

In order to make the treatment of the limiting section explicit,
we first record the elementary continuation fact that is needed. We regard $\C$ as $\R^2$;
all ODEs below are understood as equations in real variables.

\begin{lemma}[Solutions in compact sets can be continued]\label{lem:continuation}
Let $U\subset\R^d\times\R$ be open and let $b:U\to\R^d$ be a $C^1$ vector field.
If a solution of $x'=b(x,t)$ satisfies $(x(t),t)\in K$ near a finite left endpoint $a$,
where $K\subset U$ is compact, then this solution can be continued to the left beyond $a$.
\end{lemma}
\begin{proof}
$C^1$ smoothness guarantees local existence and uniqueness. Since $b$ is bounded on $K$,
there is $M<\infty$ such that $\abs{x(t)-x(u)}\le M\abs{t-u}$.
Hence $x(t)$ has a limit $x_a$ as $t\downarrow a$, and $(x_a,a)\in K\subset U$.
A local solution starting from $(x_a,a)$ exists.
After extending the original solution continuously by setting $x(a)=x_a$,
it satisfies the integral equation $x(t)=x_a+\int_a^t b(x(u),u)\,du$;
local uniqueness makes it agree with the new solution on the overlap,
so the two can be concatenated.
\end{proof}

\section{Proof of the main theorem}\label{sec:mainproof}

First let $1<q<p<2$. Necessity follows from Proposition~\ref{prop:local}. Fix $\omega_1\in A_{p,q}$;
we prove sufficiency.

If $\omega_1\in A_{p,p}$, then Theorem~\ref{thm:diagonal}
and the monotonicity of the norms on a probability space give
$$
\norm{T_{\omega_1}f}_q\le\norm{T_{\omega_1}f}_p\le\norm f_p.
$$
Henceforth assume $\omega_1\notin A_{p,p}$.

Define
$$
s_1=\sup\{s\in[q,p]:F(p,s,\omega_1)\ge0\}.
$$
By continuity of $F$ and its monotonicity in $s$ we have $q\le s_1<p$. If $F(p,s_1,\omega_1)>0$,
then by continuity one could take $s$ slightly larger and still have $F\ge0$,
contradicting the definition of the supremum. Therefore
$$
F(p,s_1,\omega_1)=0,\qquad t_1=p-s_1\in(0,p-q].
$$
This step places the target parameter, which may lie in the interior of $A_{p,q}$,
on the boundary of $A_{p,s}$ for some intermediate exponent.
Thus it suffices to prove contractivity in $L^{s_1}$ and then lower the output exponent to $q$.

In the system \eqref{eq:system} set $s(t)=p-t$ and consider the terminal value problem
\begin{equation}\label{eq:ODE}
\omega'(t)=\omega(t)\eta_*(\omega(t),p-t),
\qquad \omega(t_1)=\omega_1.
\end{equation}
The vector field is defined on the open set
$$
U=\left\{(\omega,t):
\begin{array}{l}
0<t<p-1,\quad \omega\ne0,\quad\Lambda(p,p-t,\omega)\ne0,\\
(p-t)c+2-(p-t)>0
\end{array}\right\}.
$$
Lemma~\ref{lem:tangent} guarantees that $(\omega_1,t_1)\in U$, and the vector field is $C^1$ on $U$.
Take the maximal solution through this terminal value.

In fact, the computations of \eqref{eq:D}
and of the extremal direction show that $\mathcal D_{p,s,\omega}(\eta_*)=0$ throughout $U$. Hence,
by the chain rule,
\begin{equation}\label{eq:invariant}
\frac{d}{dt}F(p,p-t,\omega(t))
=\mathcal D_{p,p-t,\omega(t)}(\eta_*)=0.
\end{equation}
Together with the terminal value this shows that the solution satisfies $F(p,p-t,\omega(t))=0$
throughout its interval of definition.

We next show that the solution exists on the whole of $(0,t_1]$.
Take any $\delta\in(0,t_1)$ and define
$$
\Sigma_\delta=
\{(\omega,t):\delta\le t\le t_1,\quad F(p,p-t,\omega)=0\}.
$$
This is a closed set; by Proposition~\ref{prop:local} we have $\abs\omega\le1$ on it,
so $\Sigma_\delta$ is compact. On $\Sigma_\delta$ we have $1<p-t<p$, and Lemma~\ref{lem:tangent}
excludes $\omega=0$, $\Lambda=0$ and $(p-t)c+2-(p-t)\le0$; therefore $\Sigma_\delta\subset U$.

If the left endpoint $t_-$ of the maximal solution were strictly larger than $\delta$,
then by \eqref{eq:invariant}
the orbit would stay entirely inside the compact set $\Sigma_\delta$ for $(t_-,t_1]$.
Lemma~\ref{lem:continuation} would then extend it to the left beyond $t_-$,
contradicting maximality. Hence $t_-\le\delta$. Since $\delta>0$ is arbitrary,
the solution exists on $(0,t_1]$, and
\begin{equation}\label{eq:boundaryflow}
\omega(t)\in\partial A_{p,p-t}\subseteq\D.
\end{equation}
Here we do not require the vector field to remain regular on the limiting section $t=0$.

Fix a non-zero two-point function $f$ and set
$$
g_t=T_{\omega(t)}f,\qquad N(t)=\norm{g_t}_{p-t}.
$$
Since $\omega(t)\ne0$, the operator $T_{\omega(t)}$ is invertible, so $g_t\ne0$.
Also $(r,z)\mapsto\abs z^r$ is jointly $C^1$ for $r>1$, so $N(t)>0$ and $N$ is a $C^1$ function.
Near a fixed $t$, \eqref{eq:ODE} gives
$$
\omega(t+\varepsilon)
=\omega(t)\bigl(1+\varepsilon\eta_*(\omega(t),p-t)+o(\varepsilon)\bigr).
$$
By the joint $C^1$ property just mentioned,
the $o(\varepsilon)$ term does not affect the derivative at $\varepsilon=0$.
Using $T_{\omega_1}T_{\omega_2}=T_{\omega_1\omega_2}$ and Lemma~\ref{lem:ITP} we obtain
\begin{equation}\label{eq:monotonic}
\frac{N'(t)}{N(t)}
=\left.\frac{d}{d\varepsilon}\right|_0
\log\norm{T_{1+\varepsilon\eta_*(\omega(t),p-t)}g_t}_{(p-t)-\varepsilon}
\le0.
\end{equation}
Hence
\begin{equation}\label{eq:compare}
\norm{T_{\omega_1}f}_{s_1}=N(t_1)\le N(t)\qquad(0<t\le t_1).
\end{equation}

By \eqref{eq:boundaryflow} and compactness of the closed unit disc,
we may choose a sequence $t_n\downarrow0$ along which $\omega(t_n)\to\omega_0\in\D$.
Since $F$ is continuous,
$$
F(p,p,\omega_0)=\lim_{n\to\infty}
F(p,p-t_n,\omega(t_n))=0.
$$
In particular $\omega_0\in A_{p,p}$. The use of Theorem~\ref{thm:diagonal}
requires only this one accumulation point;
we do not need to prove that the whole orbit has a limit at $t=0$.
By continuity of the two-point norms in $(\omega,r)$,
$$
\lim_{n\to\infty}N(t_n)=\norm{T_{\omega_0}f}_p\le\norm f_p.
$$
Taking this subsequence in \eqref{eq:compare} and using $q\le s_1$, we get
$$
\norm{T_{\omega_1}f}_q\le
\norm{T_{\omega_1}f}_{s_1}\le\norm f_p.
$$
This proves sufficiency for $1<q<p<2$.

If $2<q<p<\infty$, then $1<p'<q'<2$. The case already proved together with Lemma~\ref{lem:duality}
gives
$$
\Omega_{p,q}=\Omega_{q',p'}=A_{q',p'}=A_{p,q}.
$$
Finally, for $(a,b)=(1,0)$ both sides are equal to $1$, so the operator norm is at least $1$.
Together with the contraction bound just proved, this shows that the operator norm is exactly $1$.
This proves Theorem~\ref{thm:main}.

\section{Further results and remarks}\label{sec:remarks}

\subsection{The case \texorpdfstring{$1<q\le2\le p<\infty$}{1<q<=2<=p<infinity}}

If $1<q\le2\le p<\infty$ and $\abs\omega\le1$, then
$$
\norm{T_\omega f}_q\le\norm{T_\omega f}_2\le\norm f_2\le\norm f_p;
$$
the middle step uses $\norm{a+\omega bx}_2^2=\abs a^2+\abs{\omega b}^2$. Conversely,
taking $f=x$ gives $\abs\omega\le1$, so in this case $\Omega_{p,q}=\D$.

\subsection{Fixed points, the contraction region and equality}

For $1<s<p<2$ the four points
\begin{equation}\label{eq:anchors}
\omega=\pm1,\qquad \omega=\pm i\sqrt{p-1}
\end{equation}
all lie on $\partial A_{p,s}$. On the real axis, substitution into \eqref{eq:F}
gives $F(p,s,\pm1)=p-s-\abs{p-s}=0$; on the imaginary axis, substituting $\omega^2=-(p-1)$ gives
$$
F(p,s,\pm i\sqrt{p-1})
=p-s(p-1)-\bigl((2-p)+(2-s)(p-1)\bigr)=0.
$$
All these points satisfy $c=1$ and $\sigma=0$ in \eqref{eq:csigma}, hence $\eta_*=0$,
and they are fixed points of the characteristic flow.
As $s\uparrow p$ the quantity $\Lambda$ vanishes at the two real points;
this is exactly why the main proof uses the vector field only for $t>0$ and uses a subsequential limit at $t=0$.

On the whole of $A_{p,q}$ the operator norm equals $1$, including in the interior.
The main theorem gives the family of equality cases $(a,b)=(a,0)$ with $a\ne0$,
but it does not classify all equality cases. For example,
for $\omega=\pm1$ every function satisfying $\abs{f(1)}=\abs{f(-1)}$ has $\norm{T_\omega f}_q=\norm f_p$.
Hence one cannot infer from the contractivity that all functions with $b\ne0$ give a strict inequality.

\subsection{The boundary surface and planar flow trajectories}\label{sec:geometry}

Fix the input exponent $p$ and consider the surface in parameter space
$$
\Sigma=\{(\omega,s)\in\C\times(1,p):F(p,s,\omega)=0\}.
$$
Its three real coordinates are $(\Rea\omega,\Ima\omega,s)$;
the third coordinate is the output exponent and not the value of the norm of any function.
Fixing $s$ gives the section $\partial A_{p,s}$. By \eqref{eq:invariant}
the vector $(\omega\eta_*(\omega,s),-1)$ is tangent to $\Sigma$.
Thus the projections of the flow trajectories onto the complex plane connect boundaries
corresponding to different exponents. A planar figure is sufficient to display these trajectories.
For any fixed $f$,
along forward time $t=p-s$ the quantity $\norm{T_{\omega(t)}f}_{p-t}$ is non-increasing.
When the flow is traced backwards, the comparison is reversed:
the norm at the target does not exceed the limiting value obtained along a backward subsequence as $s$ approaches $p$.

The tangency condition for the moving boundary is
$$
\nabla_\omega F\cdot\omega'=\partial_sF.
$$
It constrains only the normal component of the planar velocity,
and it does not require the velocity to be parallel to the normal of the boundary.
Formula \eqref{eq:eta} must in addition make the control direction belong to $\mathcal K_s$,
so as to guarantee the monotonicity of the norm of every $f$. In particular,
if $\eta_*=-\alpha_*+i\beta_*$, then
\begin{equation}\label{eq:modulusphase}
\frac{d}{dt}\log\abs{\omega(t)}=-\alpha_*\le0,
\qquad \frac{d}{dt}\arg\omega(t)=\beta_*.
\end{equation}
The region $A_{p,s}$ expands as $s$ decreases,
but this does not mean that the modulus must increase along a flow trajectory:
a trajectory may start from a point of larger modulus on the original boundary,
decrease its modulus while changing its argument, and reach the boundary of the larger region.
The main proof constructs the regular flow only for $s<p$,
and as $s\uparrow p$ it uses a limit point obtained along a subsequence;
it does not require proving that the trajectory has a unique endpoint.

To compute the boundaries shown in the figure,
fix $1<q<p<2$ and write $\omega=\sqrt\rho\,e^{i\theta}$.
The determinant of the matrix in Proposition~\ref{prop:local} vanishes if and only if
$$
(q-1)\rho^2-B_{p,q}(\theta)\rho+(p-1)=0,
$$
where
$$
B_{p,q}(\theta)=\bigl((p-1)(q-1)+1\bigr)\sin^2\theta
+(p+q-2)\cos^2\theta.
$$
Taking the smaller root along the ray gives
\begin{equation}\label{eq:radial}
\rho=\frac{2(p-1)}{B_{p,q}(\theta)+
\sqrt{B_{p,q}(\theta)^2-4(p-1)(q-1)}}.
\end{equation}
$B_{p,q}$ is a convex combination of $p+q-2$ and $1+(p-1)(q-1)$,
both of which are larger than $2\sqrt{(p-1)(q-1)}$,
so the quantity under the square root is strictly positive.
For $q=p$ the quantity under the square root is non-negative and the denominator is positive,
so \eqref{eq:radial} remains valid as it stands.

\begin{figure}[p]
\centering
\includegraphics[width=.79\textwidth]{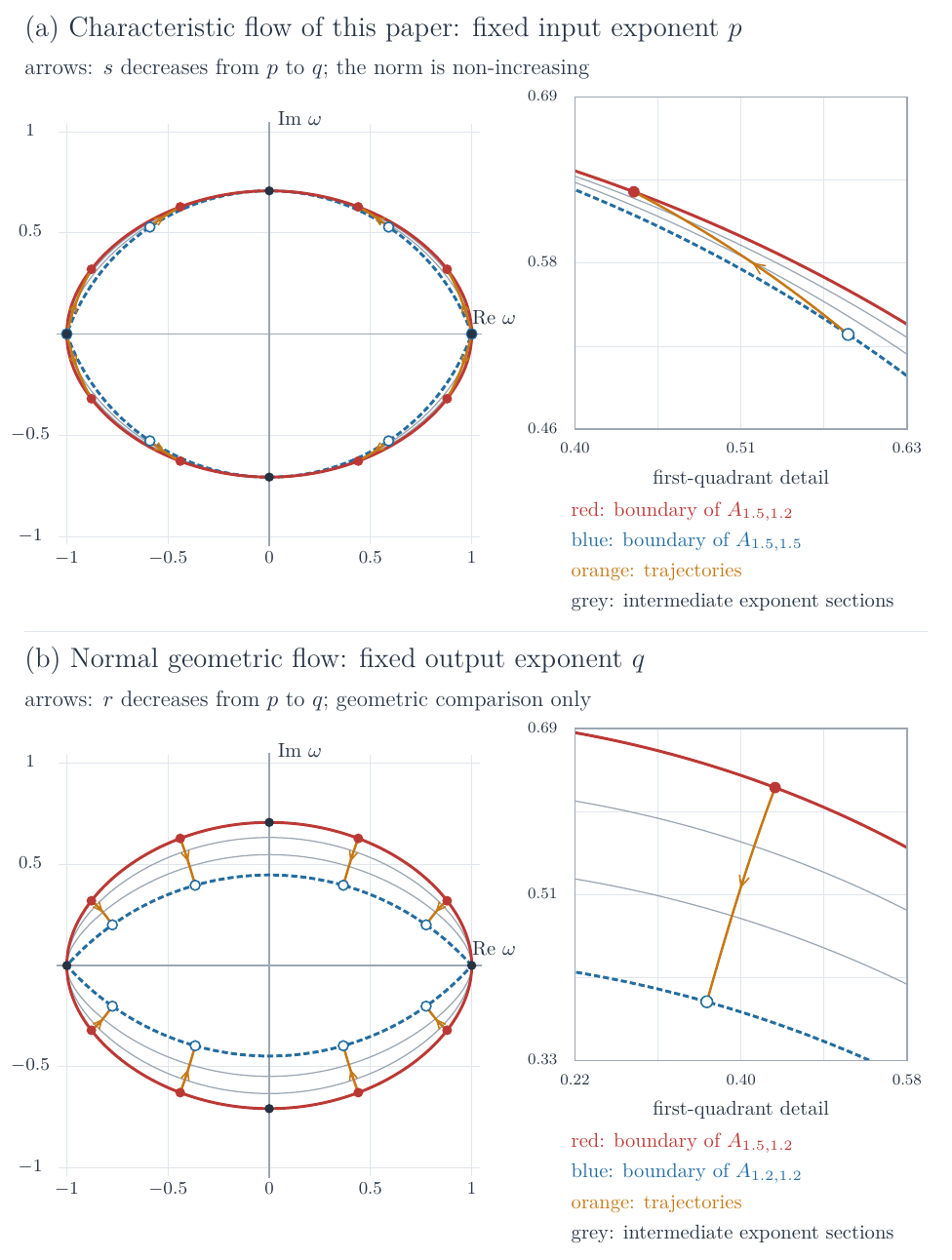}
\caption{Two types of flow trajectories for the same exponents $p=3/2$ and $q=6/5$.
(a) The characteristic flow of this paper, with fixed input exponent $p$; the
arrows point along increasing $t$ and decreasing $s$, from near
$\partial A_{p,p}$ towards $\partial A_{p,q}$.
(b) The normal geometric flow, with fixed output exponent $q$; the arrows point
along increasing $\tau$ and decreasing $r$, from $\partial A_{p,q}$ towards
$\partial A_{q,q}$. This panel is a geometric comparison only.
The solid red points are the common target points for the two flows; the open
blue points are the other endpoints of the truncated integrations. Black points
mark the target boundary points on the coordinate axes. Each right-hand inset
magnifies part of a trajectory in the first quadrant, with the same scale on
the horizontal and vertical axes.
(a) and (b) are truncated at $t=10^{-5}$ and $r-q=10^{-5}$ respectively; the
figure is not used in the rigorous proof.}\label{fig:flow}
\end{figure}

Figure~\ref{fig:flow}(a) shows the characteristic flow of this paper. Figure~\ref{fig:flow}(b)
shows a different geometric comparison: fix the output exponent $q$,
let the input exponent $r$ decrease from $p$ to $q$,
and join $\partial A_{p,q}$ to $\partial A_{q,q}$. Write $G(r,\omega)=F(r,q,\omega)$ and choose
\begin{equation}\label{eq:normalflow}
\frac{dr}{d\tau}=-1,\qquad
\frac{d\omega}{d\tau}=
\frac{\partial_rG}{\abs{\nabla_\omega G}^2}\nabla_\omega G.
\end{equation}
Then $dG/d\tau=0$, and the planar velocity is along the normal of the section boundary.
Here the planar gradient is identified with a complex number in the natural way;
if $\Lambda=(r-2)-(q-2)\omega^2$, then
$$
\partial_rG=1-\frac{\Rea\Lambda}{\abs\Lambda},\qquad
\nabla_\omega G=-2q\omega-2(2-q)\frac{\Lambda\overline\omega}{\abs\Lambda}.
$$
When $q<r<2$ and $G=0$ we have $\Lambda\ne0$ and $\omega\ne0$, and
$$
\Rea\bigl(\overline\omega\nabla_\omega G\bigr)
=-2\abs\omega^2\bigl(q+(2-q)c\bigr)<0,
$$
where $c$ is defined by \eqref{eq:csigma} at the parameters $(r,q,\omega)$.
Hence the gradient does not vanish, and \eqref{eq:normalflow}
is well defined at these boundary points. This geometric flow changes the input exponent;
it is not used in the proof of this paper and does not provide the norm-propagation conclusion needed here.
In particular,
one cannot reverse the arrows and read it as a flow deducing the target contraction from $A_{q,q}$:
along the positive imaginary axis, $\omega=i\sqrt{r-1}$,
and for $f=x$ we have $\norm{T_\omega f}_q/\norm f_r=\sqrt{r-1}$,
and this ratio increases when the flow is reversed.

\subsection{Scope of the method}\label{sec:scope}

In the boundary computation of Lemma~\ref{lem:tangent} the decisive inequality is
$$
(2-s)p-s(2-p)\cos2\theta\ge2(p-s)>0.
$$
If instead $p<s<2$, the right-hand side is no longer positive,
and $sc+2-s>0$ cannot be guaranteed near the real axis either. Hence the argument of this paper,
whose direction is $s'=-1$ and which goes back to the boundary $s=p$,
does not give a proof for $p<q$;
for the flow construction in that range see \cite{HuangIvanisvili2026}.

Nor can one dispense with Theorem~\ref{thm:diagonal}
merely by replacing the known initial set by $A_{p,2}$. For fixed $1<p<2$,
$$
A_{p,2}=\{\omega:\abs\omega\le\sqrt{p-1}\}.
$$
By real two-point hypercontractivity and the invariance of the $L^2$ norm under the phase of the multiplier,
this disc indeed gives a known contraction region
\cite{Beckner1975}.
If, however, a path makes $\norm{T_{\omega(t)}f}_{s(t)}$ non-increasing for every $f$,
then taking $f=x$ gives
$$
\norm{T_{\omega(t)}x}_{s(t)}=\abs{\omega(t)},
$$
so $\abs{\omega(t)}$ must be non-increasing.
Hence any such path starting from $A_{p,2}$ remains in that disc throughout its interval of
definition and cannot reach the points of $A_{p,q}$ whose modulus exceeds $\sqrt{p-1}$,
in particular $\omega=\pm1$. The obstruction here is the reachable set of a monotone flow,
and not merely the absence of a regular formula at $s=2$.

The real-axis boundary exhibits this conflict further. For $p<s\le2$ its positive endpoint is
$$
x_*(s)=\sqrt{\frac{p-1}{s-1}}.
$$
If $s(t)=2-t$, then
$$
\frac{d}{dt}x_*(s(t))=\frac{x_*(s(t))}{2(s(t)-1)}>0.
$$
Following this boundary requires the modulus to increase,
in contradiction with the necessary monotonicity above;
the corresponding $\eta=1/[2(s-1)]>0$ does not belong to $\mathcal K_s$ either. Thus,
under the mechanism used here, in which these norms are non-increasing for every $f$,
propagation from $A_{p,2}$ alone cannot cover the target set.
This does not exclude constructing a different monotone quantity or using other methods of proof.

As for tensorization, we make no assertion for $p>q$ analogous to the one for $p\le q$.
The mixed-norm exchange inequalities used in induction on the coordinates have a different direction when $p>q$.
The continuous Gaussian problem also has different boundedness properties at the boundary:
the criterion for $q<p$ in Weissler \cite[Theorem~1]{Weissler1979}
requires a strict determinant condition when $\Rea(\omega/(1-\omega^2))\ne0$, so, in general,
non-axial Weissler boundary points do not satisfy the Gaussian boundedness condition.
Thus two-point boundary contraction cannot be transferred to Gaussian boundary contraction without further argument.

\section*{Acknowledgements}

The author acknowledges the use of AI-assisted tools in the preparation of this manuscript.
All mathematical results have been verified by the author,
who takes full responsibility for the content of the paper.

\end{document}